\documentclass[11pt, a4paper]{amsart}
\usepackage[utf8]{inputenc}
\usepackage{amsfonts, amssymb, amsmath, textcomp, amsthm, xcolor}

\def\R{\mathbb{R}}

\def\nint{\mathop{\diagup\kern-13.0pt\int}}
\def\Z{\mathbb{Z}}
\def\T{\mathbb{T}}

\def\beq{\begin{equation}}
\def\endeq{\end{equation}}
\def\bg{\begin{gathered}}
\def\eg{\end{gathered}}

\newtheorem{thm}{Theorem}[section]

\newtheorem{cor}[thm]{Corollary}

\newtheorem{rem}[thm]{Remark}

\title[On a free Schr\"{o}dinger solution studied by BBCRV]{On a free Schr\"{o}dinger solution studied by Barcel\'{o}--Bennett--Carbery--Ruiz--Vilela}
\author{Xiumin Du, Yumeng Ou, Hong Wang and Ruixiang Zhang}
\date{\today}
\thanks{XD is supported by NSF DMS-2107729 (transferred from DMS-1856475) and Sloan Research Fellowship. YO is supported by NSF DMS-2142221 and NSF DMS-2055008. HW is supported by NSF DMS-2141426. RZ is supported by NSF DMS-2207281 (transferred from DMS-1856541).}

\begin{document}

\maketitle

\begin{abstract}
We present a free Schr\"{o}dinger solution studied by Barcel\'{o}--Bennett--Carbery--Ruiz--Vilela and show why it can be viewed as a sharp example for the recently discovered refined decoupling theorem.
\end{abstract}

\section{Introduction}

Let $d\geq 2$. Consider the free Schr\"{o}dinger equation:

\begin{equation}\label{schrodingereq}
\left\{
\begin{array}{l}
\mathrm{i}u_{t} - \Delta_x (u) = 0,\\
u(x, 0) = f(x)
\end{array} \right.
\end{equation} for $(x, t) \in \R^{d-1} \times \R$. By taking the Fourier transform of both sides, we know $\text{supp}\,\hat{u} \subset \widetilde{P^{d-1}}$. Here $\widetilde{P^{d-1}}$ is the paraboloid:
$$\widetilde{P^{d-1}} = \{\xi_d = |\xi'|^2\},\quad \xi=(\xi',\xi_d):=(\xi_1,\ldots, \xi_{d-1}, \xi_d)\in \mathbb{R}^d.$$

Because of the above property, such functions $u$ are closely related to the \emph{Fourier restriction theory} and have been extensively studied by Fourier analysts. In light of the Littlewood-Paley decomposition, people are often interested in functions $g$ on $\R^n$ such that $\text{supp}\, \hat{g}$ is in the \emph{truncated paraboloid} $$P^{d-1} = \{\xi_d = |\xi'|^2,\, |\xi_j| \leq 1, \forall 1 \leq j \leq d-1\}.$$ Fourier analysts are then interested in $L^p \to L^q$ estimates of such  $g$ on certain subsets of $\R^d$, and a great amount of related recent progress has been made.

In this note, we first review an example of such a function $g$ studied by Barcel\'{o}--Bennett--Carbery--Ruiz--Vilela in \cite{barcelo2007some}. Next, we present a recent result known as ``refined decoupling'' (proved independently by Guth--Iosevich--Ou--Wang \cite{guth2020falconer} and Du--Zhang). Refined decoupling has seen powerful applications in recent years such as in the Falconer distance problem \cite{guth2020falconer, du2021improved} and small cap decouplings \cite{DGW}. It would thus be interesting to know various sharp examples for this estimate. In this note, we show that  Barcel\'{o}-Bennett--Carbery--Ruiz--Vilela's free Schr\"{o}dinger solutions are always (almost) sharp examples for refined decoupling.
 Moreover, we show in the end that this example is also sharp for an $L^2$ estimate used as a key step in many recent arguments studying the Falconer distance problem.

 \begin{rem}\label{rem1}
 Historically, \cite{barcelo2007some} introduced this example and generalizations to provide useful test cases for
$L^2$-average decay estimates of Fourier transforms of
fractal measures. Therefore, it is perhaps not surprising that the example may be relevant for testing against other related results such as refined decoupling.
\end{rem}
 
 \begin{rem}\label{rem2}  In  \cite{CIW}, inequality (4) is essentially a restatement of refined decoupling and it was remarked (Remark 1.3) that Knapp examples make refined decoupling (almost) sharp too.
To put \cite{CIW} into historical context, Guth \cite{Gu22} provided another different example  capturing limits of decoupling. Guth's example was worked out carefully in \cite{CIW} to show sharpness of their study of Mizohata-Takeuchi conjecture using refined decoupling estimate. 
\end{rem}

\section{Barcel\'{o}--Bennett--Carbery--Ruiz--Vilela's free Schr\"{o}dinger solution}\label{BBCRVgsec}

Let $0<\sigma < 1/2$ and $R>1$ be fixed parameters. Let $\mathrm{d}\omega$ be the hypersurface measure on $P^{d-1}$. Barcel\'{o}--Bennett--Carbery--Ruiz--Vilela's free Schr\"{o}dinger solution is a function $g$ such that $$\hat{g} (\xi) = h(\xi)\, \mathrm{d}\omega$$ where $$h(\xi) = \sum_{\substack{l_1,\ldots, l_{d-1} \in \Bbb{Z},\\ 1 \leq l_1,\ldots, l_{d-1} <R^{\sigma}}} 1_{(l_1 R^{-\sigma} -R^{-1}, l_{1}R^{-\sigma} + R^{-1}) \times \cdots \times (l_{d-1}R^{-\sigma} -R^{-1}, l_{d-1}R^{-\sigma} + R^{-1})}.$$ 

We now state the most relevant properties of the above function $g$ here. By elementary computations, one can check that $|g|\sim R^{(d-1)(\sigma -1)}$ at all points of the form $(n_1 R^{\sigma}, n_2 R^{\sigma}, \ldots, n_{d-1} R^{\sigma}, n_d R^{2\sigma})$ inside the ball $B_{c_d R}$ of radius $c_d R$, where $n_j\in \mathbb{Z}$ and $c_d>0$ is a small constant only depending on the ambient dimension. Note that $R^{(d-1)(\sigma -1)}$ is comparable to $\|g\|_{\infty}$ by triangle inequality. Moreover, $|g|\sim R^{(d-1)(\sigma -1)}$ inside a ball of radius $\sim_d 1$ around every point above by a similar computation. We refer the reader to \cite{barcelo2007some} for more detailed justification of these facts.

\section{Wave packet decomposition}

In order to introduce the refined decoupling inequality, we first briefly recall the \emph{wave packet decomposition}, a standard tool for analyzing functions with Fourier support in $P^{d-1}$. Here, we present (the rescaled version of) the wave packet decomposition used in \cite{guth2020falconer}.

Fix a parameter $R>1$. Decompose $P^{d-1}$ into pieces $\theta$ such that the projection of each $\theta$ onto the hyperplane of the first $d-1$ coordinates is a square of side length $R^{-\frac{1}{2}}$. By elementary differential geometry, each $\theta$ is contained in a box of dimensions $R^{-\frac{1}{2}} \times \cdots \times R^{-\frac{1}{2}} \times R^{-1}$. Pick such a box and let $T_{\theta}$ be its dual box centered at the origin. Note that $T_{\theta}$ is roughly a tube of thickness $R^{\frac{1}{2}}$ and length $R$. Let $B_R \subset \R^d$ be the ball centered at the origin with radius $R$. Tile $B_R$ by translations of $T_{\theta}$ and call this family $\T_{\theta}$.

For a function $v$ whose Fourier support is in $P^{d-1} \subset \R^d$, one can decompose $$v = \sum_{\theta, T:\, T \in \T_{\theta}} v_{\theta, T}$$ inside $B_R$ such that:

\begin{itemize}
\item Each $\widehat{v_{\theta, T}}$ is supported in a box of dimensions $\sim R^{-\frac{1}{2}} \times \cdots \times R^{-\frac{1}{2}} \times R^{-1}$ containing $\theta$.
\item Each $v_{\theta, T}$ (known as a wave packet) is morally supported in $T$ and rapidly decays outside of it.
\item Each $|v_{\theta, T}|$ is morally a constant on $T$ and we will call this constant the \emph{magnitude} of $v_{\theta, T}$.
\item Different $v_{\theta, T}$ are morally $L^2$-orthogonal on every $R^{\frac{1}{2}}$-ball. This property is known as \emph{local orthogonality} and is very useful in Fourier restriction type problems. We do not need its detailed description here.
\end{itemize}

For each $\theta$, we also define  $$v_{\theta} = \sum_{T \in \T_{\theta}} v_{\theta, T}.$$

\section{Refined decoupling and its sharpness}

\subsection{The refined decoupling theorem}

We can now state the refined decoupling theorem:

\begin{thm}[Refined decoupling \cite{guth2020falconer}]\label{refdecthm}
Suppose $\text{supp}\, \hat{v} \subset P^{d-1}$ and $R>1$. Suppose that in the wave packet decomposition of $v$ in $B_R$, every two wave packets have comparable magnitudes. Let $X \subset B_R$ such that each $x \in X$ hits the essential support of $\leq M$ wave packets, then for $p = \frac{2(d+1)}{d-1}$,
\begin{equation}\label{refdecineq}
\|v\|_{L^p (X)} \lesssim_{\varepsilon} R^{\varepsilon} M^{\frac{1}{2} - \frac{1}{p}} \left(\sum_{\theta} \|v_{\theta}\|_{L^p (w_{B_R})}^p\right)^{\frac{1}{p}}.
\end{equation}
\end{thm}

Here we have a  weight $w_{B_R}$ included on the right hand side for technical reasons. It behaves like $1_{B_R}$ but has a rapidly decaying tail. Morally one can think of $\|\cdot\|_{L^p (w_{B_R})}$ as $\|\cdot\|_{L^p (B_R)}$. Theorem \ref{refdecthm} is named refined decoupling because it is a refinement of the celebrated Bourgain--Demeter decoupling theorem for paraboloids \cite{bourgain2015proof}.

We remark that the assumption that all wave packets have comparable magnitudes in the theorem is usually harmless. In applications, one can usually reduce a general situation to this case by dyadic pigeonholing. 

In \cite{guth2020falconer}, Theorem \ref{refdecthm} is one of the central ingredients the authors use to make progress on the Falconer distance conjecture in $\R^2$. The theorem is also useful in other problems of similar flavors such as Schr\"{o}dinger maximal function estimates.

\subsection{An almost sharp example for Theorem \ref{refdecthm}} The function $g$ we discussed in \S \ref{BBCRVgsec} is an (almost) sharp example for Theorem \ref{refdecthm}, as we explain below.

Note that for the function $g$, if we take $X$ to be the $1$-neighborhood of $\{(n_1 R^{\sigma}, n_2 R^{\sigma}, \ldots, n_{d-1} R^{\sigma}, n_d R^{2\sigma}): n_1, \ldots, n_d \in \Z\} \bigcap B_R^d$, then since the Fourier support of $g$ only intersects $\sim R^{(d-1)\sigma}$ $\theta$'s and that the supports of wave packets from one $\theta$ are essentially non-overlapping,  the relevant $M$ is $\lesssim R^{(d-1)\sigma}$. In fact, since $g$ attains almost its maximal possible value at each $(n_1 R^{\sigma}, n_2 R^{\sigma}, \ldots, n_{d-1} R^{\sigma}, n_d R^{2\sigma})$, one can further see that $M$ is indeed $\sim R^{(d-1)\sigma}$, but even without this stronger observation one can still see sharpness from the computation below.

Recall that $p = \frac{2(d+1)}{d-1}$. Because of the property that $g$ is almost the largest possible at each $(n_1 R^{\sigma}, n_2 R^{\sigma}, \ldots, n_{d-1} R^{\sigma}, n_d R^{2\sigma})$, as discussed in  \S \ref{BBCRVgsec}, we see that 
\begin{equation}\label{gleftrefdec}
\|g\|_{L^p (X)} \sim R^{(d-1)(\sigma - 1)} |X|^{1/p} \sim R^{(d-1)(\sigma - 1)} R^{\frac{d-(d+1)\sigma}{p}} \sim R^{\frac{d-1}{2}\sigma -\frac{(d-1)(d+2)}{2(d+1)}}.
\end{equation}

Let us look at the right hand side of (\ref{refdecineq}). Each $g_{\theta}$ is supported in a ball of radius $\sim R^{-1}$, so one can see that each $|g_{\theta}|$ is $\sim R^{-(d-1)}$ on a ball of radius $\sim R$ centered at the origin and has the same value on the whole space as an upper bound. Hence each $$\|g_{\theta}\|_{L^p (w_{B_R})} \sim R^{- (d-1)+ \frac{d}{p}} \sim R^{-\frac{(d-1)(d+2)}{2(d+1)}}$$ and the right hand side  of (\ref{refdecineq}) for $v=g$ is

\begin{equation}\label{grightrefdec}
\sim R^{\varepsilon} \cdot R^{\frac{(d-1)\sigma}{d+1}}\cdot R^{\frac{(d-1)\sigma}{p}} \cdot R^{-\frac{(d-1)(d+2)}{2(d+1)}} \sim R^{\frac{d-1}{2}\sigma -\frac{(d-1)(d+2)}{2(d+1)}+\varepsilon}.
\end{equation}

We see that the right hand sides of (\ref{gleftrefdec}) and (\ref{grightrefdec}) match except for the $R^\varepsilon$-loss, showing that $g$ is an almost sharp example for Theorem \ref{refdecthm}.

Since the decoupling theorem for the paraboloid by Bourgain--Demeter is weaker than Theorem \ref{refdecthm}, we see that this function $g$ is also sharp for Bourgain--Demeter's decoupling  theorem.

As mentioned above, the refined decoupling Theorem \ref{refdecthm} has other consequences that are useful in problems in geometric measure theory such as the Falconer distance problem. For example, by dyadic pigeonholing and H\"{o}lder's ineuqality, it implies the following corollary:

\begin{cor}\label{L2estimatecor}
Let $R>1$ and $\alpha >0$. Suppose a set $Y \subset B_R$ is a union of lattice 1-cubes with the following ``fractal structure (between scale $1$ and $R$)'': $$|B_r \bigcap Y| \lesssim r^{\alpha},\quad \forall 1 \leq r \leq R,\, \forall B_r \subset B_{R}.$$ Suppose $\text{supp}\, \hat{v} \subset P^{d-1}$ with $\hat{v} = \varphi \,\mathrm{d}\omega$ such that in the wave packet decomposition of $v$ in $B_R$, each wave packet has its essential support hits a $\lesssim R^{-\frac{(d-1)}{2}}$ fraction of unit cubes in $Y$. Then

\begin{equation}\label{L2estimateineq}
\|v\|_{L^2 (Y)} \lesssim_{\varepsilon} R^{\frac{1}{d+1} (\alpha -\frac{d-1}{2})+\varepsilon} \|\varphi\|_{L^2 (
\mathrm{d}\omega)}.
\end{equation}
\end{cor}

By a similar computation, one can see that if we take $\alpha = d- (d+1) \sigma$ and $v$ to be the function $g$ in \S \ref{BBCRVgsec}, one has again an almost sharp example for Corollary \ref{L2estimatecor}. In this example, both sides of (\ref{L2estimateineq}) are comparable to or close to $R^{\frac{d-3}{2}\sigma -\frac{d-2}{2}}$.

Theorem \ref{refdecthm} and various versions of Corollary \ref{L2estimatecor} were used in \cite{guth2020falconer} and later works such as \cite{du2021improved} to make progress towards the Falconer distance conjecture. The sharpness of these two propositions showed in this section suggests that to make further progress beyond e.g. \cite{guth2020falconer}, one  either  has to sharpen other components in these two papers or to design new approaches.

\bibliography{ref}{}
\bibliographystyle{alpha}
\vspace{1cm}
	
	\noindent Xiumin Du. Northwestern University.
	Email address: xdu@northwestern.edu\\

	\noindent Yumeng Ou. University of Pennsylvania.
	Email address: yumengou@sas.upenn.edu\\

\noindent Hong Wang. UCLA. Email address: hongwang@math.ucla.edu\\

\noindent Ruixiang Zhang.  UC Berkeley. Email address: ruixiang@berkeley.edu

\end{document}